\documentclass[12pt]{article}
\usepackage{amsmath}
\usepackage{epsf}
\usepackage{amssymb}
\usepackage{graphicx}
\usepackage{ifpdf}
\usepackage{caption}

\usepackage{array}
\usepackage{booktabs} 
\usepackage{multirow}

\newtheorem{theorem}{Theorem}[section]

\newtheorem{definition}{Definition}[section]
\newtheorem{lemma}{Lemma}[section]
\newtheorem{corollary}{Corollary}[section]

\begin{document}
\renewcommand{\thefootnote}{\fnsymbol{footnote}}

\begin{center}
{\bf \Large The generalized 3-connectivity of a family of regular networks}\footnote{This
work was supported by National Science Foundation of China(No.12271157), Natural Science Foundation of Hunan Province (No.2022JJ30028) and Hunan Education Department Foundation(No.21C0762).}
\vskip 5mm

{{\bf Jing Wang$^1$, Xidao Luan$^2$, Yuanqiu Huang$^3$}\\[2mm]
$^1$ School of Mathematics, Changsha University, Changsha 410022, China\\
$^2$ School of Computer Technology and Science, Changsha University, Changsha 410022, China\\
$^3$ School of Mathematics, Hunan Normal University, Changsha 410081, China}\\[6mm]
\end{center}
\date{}

\noindent{\bf Abstract}\; The generalized $k$-connectivity of a graph $G$, denoted by $\kappa_k(G)$, is the minimum number of internally edge disjoint $S$-trees for any $S\subseteq V(G)$ with $|S|=k$. The generalized $k$-connectivity is a natural extension of the classical connectivity and plays a key role in applications related to the modern interconnection networks. In this paper, we firstly introduce a family of regular networks $H_n$ that can be obtained from several subgraphs $G_n^1, G_n^2, \cdots, G_n^{t_n}$ by adding a matching, where each subgraph $G_n^i$ is isomorphic to a particular graph $G_n$ ($1\le i\le t_n$). Then we determine the generalized 3-connectivity of $H_n$. As applications of the main result, the generalized 3-connectivity of some two-level interconnection networks, such as the hierarchical star graph $HS_n$, the hierarchical cubic network $HCN_n$ and the hierarchical folded hypercube $HFQ_n$, are determined directly.

{\bf MR(2000) Subject Classification} 05C40, 05C05

\section{Introduction}
\label{secintro}

With rapid development and advances of very large scale integration technology and wafer-scale integration technology, multiprocessor systems have been widely designed and used in our daily life. It is well known that the underlying topology of the multiprocessor systems can be modelled by a connected graph $G=(V(G),E(G))$, where $V(G)$ is the set of processors and $E(G)$ is the set of communication links of multiprocessor systems.

A subset $S\subseteq V(G)$ of a connected graph $G$ is called a {\it vertex-cut} if $G\setminus S$ is disconnected or trivial. The {\it connectivity} $\kappa(G)$ of $G$ is defined as the minimum cardinality over all vertex-cuts of $G$. The connectivity $\kappa(G)$ of $G$ is an important measurements for fault tolerance of the network since the larger $\kappa(G)$ is, the more reliable the network is. A well known theorem of Whitney \cite{Whitney1932} provides an equivalent definition of connectivity. For each 2-subset $S=\{x,y\}\subseteq V(G)$, let $\kappa(S)$ denote the maximum number of internally disjoint ($x,y$)-paths in $G$. Then
\begin{equation*}
\kappa(G)=\min\{\kappa(S) |S\subseteq V(G)\; {\rm and} \; |S|=2\}.
\end{equation*}

The generalized $k$-connectivity, which was introduced by Chartrand et al. \cite{Chartrand1984}, is a strengthening of connectivity and can be served as an essential parameter for measuring reliability and fault tolerance of the network. Let $G=(V(G),E(G))$ be a simple graph, $S$ be a subset of $V(G)$. A tree $T$ in $G$ is called an $S$-{\it tree}, if $S\subseteq V(T)$. The trees $T_1, T_2, \cdots, T_r$ are called {\it internally edge disjoint $S$-trees} if $V(T_i)\cap V(T_j)=S$ and $E(T_i)\cap E(T_j)=\emptyset$ for any integers $1\le i\ne j\le r$. $\kappa_G(S)$ denote the maximum number of internally edge disjoint $S$-trees. For an integer $k$ with $2\le k\le |V(G)|$, the {\it generalized $k$-connectivity} of $G$, denoted by $\kappa_k(G)$, is defined as
\begin{equation*}
\kappa_k(G)=\min\{\kappa_G(S) |S\subseteq V(G)\; {\rm and} \; |S|=k\}.
\end{equation*}

The generalized 2-connectivity is exactly the traditional connectivity. Over the past few years, research on the generalized connectivity has received meaningful progress. Li et al. \cite{SLi2012n} derived that it is NP-complete for a general graph $G$ to decide whether there are $l$ internally edge disjoint trees connecting $S$, where $l$ is a fixed integer and $S\subseteq V(G)$. Authors in \cite{HZLi2014,SLi2010} investigated the upper and lower bounds of the generalized connectivity of a general graph $G$.

Many authors tried to study exact values of the generalized connectivity of graphs. The generalized $k$-connectivity of the complete graph, $\kappa_k(K_n)$, was determined in \cite{Chartrand2010} for every pair $k,n$ of integers with $2\le k\le n$. The generalized $k$-connectivity of the complete bipartite graphs $K_{a,b}$ are obtained in \cite{SLi2012b} for all $2\le k\le a+b$. Apart from these two results, the generalized $k$-connectivity of other important classes of graphs, such as, Cartesian product graphs \cite{HZLi2012,HZLi2017}, hypercubes \cite{HZLi2012,SLin2017}, several variations of hypercubes \cite{ZhaoHao20191,ZhaoHao20192,Wei2021,Wang2021}, Cayley graphs \cite{SLi2017,ZHao20193,SLi2016,Hao20191}, have draw many scholars' attention. So far, as we can see, the results on the generalized $k$-connectivity of network are almost about $k=3$.

For large systems, it is desirable to have a cluster-based or hierarchical interconnection network, in which lower level networks support local communication, and higher level networks support remote communication. The hierarchical star graph $HS_n$\cite{HSn2005}, the hierarchical cubic network $HCN_n$ \cite{HCNn1995} and the hierarchical folded hypercube $HFQ_n$\cite{HFQn2000} are three kinds of two-level interconnection networks.  All of them are regular and have been used to design various commercial multiprocessor machines since they possess many desirable properties, such as low degree, small diameter, an so on.

The paper is organized as follows. Section \ref{secpreli} gives some necessary preliminaries. In Section \ref{secHGn}, we firstly introduce a family of regular networks $H_n$ that can be obtained from several subgraphs $G_n^1, G_n^2,\cdots, G_n^{t_n}$, where each subgraph $G_n^i$ is isomorphic to a particular graph $G_n$ ($1\le i\le t_n$). The generalized 3-connectivity of $H_n$ is then studied.  As applications of the main result, the generalized 3-connectivity of three two-level interconnection networks, such as the hierarchical star graph $HS_n$, the hierarchical cubic network $HCN_n$ and the hierarchical folded hypercube $HFQ_n$, are determined in Section \ref{secappli}. In Section \ref{secconclusion}, the paper is concluded.

\section{Preliminaries}\label{secpreli}

This section is dedicated to introduce some necessary preliminaries. We only consider a simple, connected graph $G=(V(G),E(G))$ with $V(G)$ be its vertex set and $E(G)$ be its edge set. For a vertex $x\in V(G)$, the {\it degree} of $x$ in $G$, denoted by ${\rm deg}_G(x)$, is the number of edges of $G$ incident with $x$. Denote $\delta(G)$ the {\it minimum degree} of vertices of $G$. We can abbreviate $\delta(G)$ to $\delta$ if there is no confusion. A graph is $d$-{\it regular} if ${\rm deg}_G(x)=d$ for every vertex $x\in V(G)$. For a vertex $x\in V(G)$, we use $N_G(x)$ to denote the neighbour vertices set of $x$ and $N_G[x]$ to denote $N_G(x)\cup \{x\}$. Let $V'\subseteq V(G)$, denote by $G\backslash V'$ the graph obtained from $G$ by deleting all the vertices in $V'$ together with their incident edges. 

Let $P$ be a path in $G$ with $x$ and $y$ be its two terminal vertices, then $P$ is called an $(x,y)$-{\it path}. Two ($x,y$)-paths $P_1$ and $P_2$ are {\it internally disjoint} if they have no common internal vertices, that is, $V(P_1)\cap V(P_2)=\{x, y\}$.

Li et al. \cite{SLi2010} gave an upper and lower bound of $\kappa_3(G)$ for a general graph $G$.

\begin{lemma}\label{lemupperK}(\cite{SLi2010})
Let $G$ be a connected graph with minimum degree $\delta$. If there are two adjacent vertices of degree $\delta$, then $\kappa_3(G)\le \delta-1$.
\end{lemma}

\begin{lemma}\label{lemkr}(\cite{SLi2010})
Let $G$ be a connected graph with $n$ vertices. For every two integers $k$ and $r$ with $k\ge 0$ and $r\in\{0,1,2,3\}$, if $\kappa(G)=4k+r$, then $\kappa_3(G)\ge 3k+\lceil \frac{r}{2}\rceil$.
\end{lemma}

\begin{lemma}\label{lemxypath} (\cite{Bondy})
Let $G$ be a $k$-connected graph, and let $x$ and $y$ be a pair of distinct vertices of $G$. Then there exist $k$ internally disjoint ($x,y$)-paths in $G$.
\end{lemma}

\begin{lemma}\label{lemKfan} (\cite{Bondy})
Let $G$ be a $k$-connected graph, let $x$ be a vertex of $G$ and let $Y\subseteq V(G)\backslash\{x\}$ be a set of at least $k$ vertices of $G$. Then there exists a $k$-fan in $G$ from $x$ to $Y$, that is, there exists a family of $k$ internally disjoint ($x,Y$)-paths whose terminal vertices are distinct in $Y$.
\end{lemma}

\section{The definition of $H_n$}\label{secHGn}

Let $[n]=\{1,2,\cdots, n\}$. Firstly, we introduce a family of regular graphs $H_n$ which can be constructed from $t_n$ different subgraphs $G_n^1, G_n^2,\cdots, G_n^{t_n}$, each of which is isomorphic to a particular graph $G_n$.

\begin{definition}\label{defHGn}
For integers $d$ and $t_n$ satifying $t_n\ge d+3$, let $G_n$ be a given $d$-regular $d$-connected graph with $t_n$ vertices, moreover, $\kappa_3(G_n)=d-1$. Set $G_n^1, G_n^2,\cdots,  G_n^{t_n}$ be $t_n$ different copies of $G_n$. Define $H_n$ be a ($d+1$)-regular graph obtained from $G_n^1\cup G_n^2\cup\cdots \cup G_n^{t_n}$ by adding  $\frac{1}{2}t_n^2$ edges satisfying the following two conditions:\\
(1) for each vertex $x\in V(G_n^i)$ ($1\le i\le t_n$), it has exactly one neighbour outside $G_n^i$, which is called the out-neighbour of $x$ and denoted by $\hat{x}$;\\
(2) for $1\le i\ne j\le t_n$, there is one or two cross edges between different subgraphs $G_n^i$ and $G_n^j$.\\
We write the construction of $H_n$ symbolically as $H_n=G_n^1\oplus G_n^2\oplus \cdots \oplus G_n^{t_n}$. Each $G_n^i$ is called a cluster of $H_n$ ($1\le i\le t_n$).
\end{definition}

\begin{lemma}\label{lemHG1}
Let $H_n=G_n^1\oplus G_n^2\oplus \cdots \oplus G_n^{t_n}$ and $H=H_n\backslash V(G_n^i)$, where $G_n^i$ is any cluster of $H_n$, $1\le i\le t_n$. Then $\kappa(H)=d$.
\end{lemma}

\noindent {\bf Proof} \; Firstly, $\kappa(H)\le \delta(H)=d$. To obtain the reverse inequality, we need to show that there are $d$ internally disjoint ($x,y$)-paths for any two vertices $x$ and $y$ in $H$. The following two cases are considered.

\vskip 2mm
{\bf Case 1.}\; Both $x$ and $y$ belong to a same cluster, say $G_n^{1}$.

By Definition \ref{defHGn}, there are $d$ internally disjoint ($x,y$)-paths in $G_n^{1}$ since $\kappa(G_n^{1})=\kappa(G_{n})=d$.

\vskip 2mm
{\bf Case 2.}\; $x$ and $y$ belong to different clusters.

W.l.o.g., assume that $H=H_n\backslash V(G_n^1)$, $x\in V(G_n^{2})$ and $y\in V(G_n^{3})$. According to Definition \ref{defHGn}, there exists an edge $u_i\hat{u}_i$ between $G_n^{2}$ and $G_n^{i+3}$, where $u_i\in V(G_n^{2})$ and $\hat{u}_i\in V(G_n^{i+3})$ for $1\le i\le d$. This is possible since $t_n\ge d+3$ by Definition \ref{defHGn}. Analogously, there is an edge $w_i\hat{w}_i$ between $G_n^{3}$ and $G_n^{i+3}$, where $w_i\in V(G_n^{3})$ and $\hat{w}_i\in V(G_n^{i+3})$ for $1\le i\le d$.

Let $U=\{u_1,\cdots, u_{d}\}$  and $W=\{w_1,\cdots, w_d\}$. It is seen that $|U|=|W|=d$. By Lemma \ref{lemKfan}, for $1\le i\le d$, there is a family of $d$ internally disjoint ($x,U$)-paths $Q_1,\cdots, Q_{d}$ in $G_n^{2}$ such that $u_i\in V(Q_i)$ and a family of $d$ internally disjoint ($y,W$)-paths $R_1,\cdots, R_{d}$ in $G_n^{3}$, where $w_i\in V(R_i)$.

For $1\le i\le d$, there is a ($\hat{u}_i, \hat{w}_i$)-path $\widetilde{P}_i$ in $G_n^{i+3}$ since $G_n^{i+3}$ is connected. Let
\begin{equation*}
P_i=Q_i\cup R_i\cup \widetilde{P}_i\cup\{u_i\hat{u}_i, w_i\hat{w}_i\}, \;\;\;\; 1\le i\le d.
\end{equation*}
Then $P_1, \cdots, P_d$ are $d$ internally disjoint ($x,y$)-paths in $H$.    \hfill$\Box$

\begin{theorem}\label{thmHGn2}
Let $H_n=G_n^1\oplus G_n^2\oplus \cdots \oplus G_n^{t_n}$. Then $\kappa_3(H_n)=d$.
\end{theorem}

\noindent {\bf Proof} \; Firstly, $\kappa_3(H_n)\le \delta(H_n)-1=d$ by Lemma \ref{lemupperK} and Definition \ref{defHGn}. Now we are going to prove the reverse inequality. Let $S=\{x,y,z\}$ be any 3-subset of $V(H_n)$.

\vskip 2mm
{\bf Case 1.}\; $x, y$ and $z$ belong to a same cluster of $H_n$, say $G_n^1$.

By Definition \ref{defHGn}, there are ($d-1$)-internally edge disjoint $S$-trees $T_1, \cdots, T_{d-1}$ in $G_n^1$ since $\kappa_3(G_n^1)=\kappa_3(G_n)=d-1$. Recall that $\hat{x}$, $\hat{y}$ and $\hat{z}$ are out-neighbours of $x,y$ and $z$, respectively. It follows from Lemma \ref{lemHG1} that there is an $\{\hat{x}, \hat{y}, \hat{z}\}$-tree $\widetilde{T}_{d}$ in $H_n\backslash V(G_n^1)$ since $H_n\backslash V(G_n^1)$ is connected. Let
\begin{equation*}
T_{d}=\widetilde{T}_{d}\cup \{x\hat{x}, y\hat{y}, z\hat{z}\}.
\end{equation*}
Then $T_1, \cdots, T_{d-1}, T_{d}$ are $d$-internally edge disjoint $S$-trees in $H_n$.

\vskip 2mm
{\bf Case 2.}\; $x, y$ and $z$ belong to two different clusters of $H_n$.

W.l.o.g., assume that $\{x,y\}\subseteq V(G_n^1)$ and $z\in V(G_n^2)$. By Definition \ref{defHGn}, there exist $d$ internally disjoint ($x,y$)-paths $P_1, \cdots, P_{d}$ in $G_n^1$ since $\kappa(G_n^1)=\kappa(G_n)=d$. Let $u_i$ ba a neighbour of $x$ with $u_i\in V(P_i)$ for $1\le i\le d$. It is possible that $y\in \{u_1, \cdots, u_{d}\}$. This possibility doesn't affect the following discussions.

Let $\hat{U}=\{\hat{u}_1,\cdots, \hat{u}_{d}\}$. Clearly, $\hat{U}\subseteq V(H_n)\backslash V(G_n^1)$ and $|\hat{U}|=d$. According to Lemma \ref{lemHG1} and Lemma \ref{lemKfan}, there is a family of $d$ internally disjoint ($z,\hat{U}$)-paths $Q_1, \cdots, Q_d$ in $H_n\backslash V(G_n^1)$ where $\hat{u}_i\in V(Q_i)$, $1\le i\le d$.

For $1\le i\le d$, let
\begin{equation*}
T_{i}=P_i\cup Q_i\cup \{u_i\hat{u}_i\}.
\end{equation*}
Then $T_1, \cdots,  T_{d}$ are $d$-internally edge disjoint $S$-trees in $H_n$.

\vskip 2mm
{\bf Case 3.}\; $x, y$ and $z$ belong to three different clusters of $H_n$.

W.l.o.g., assume that $x\in V(G_n^1)$, $y\in V(G_n^2)$ and $z\in V(G_n^3)$. By Definition \ref{defHGn}, there is an edge $u_i\hat{u}_i$ between subgraphs $G_n^1$ and $G_n^{i+3}$ where $u_i\in V(G_n^1)$ and $\hat{u}_i\in V(G_n^{i+3})$, $1\le i \le d$. This is possible since $t_n\ge d+3$. Similarly, for $1\le i\le d$, there is an edge $v_i\hat{v}_i$ between subgraphs $G_n^2$ and $G_n^{i+3}$ where $v_i\in V(G_n^2)$ and $\hat{v}_i\in V(G_n^{i+3})$, there is an edge $w_i\hat{w}_i$ between subgraphs $G_n^3$ and $G_n^{i+3}$ where $w_i\in V(G_n^3)$ and $\hat{w}_i\in V(G_n^{i+3})$.

Combined with Definition \ref{defHGn} and Lemma \ref{lemKfan}, there is a $d$-fan $P_1, \cdots, P_d$ in $G_n^1$ from $x$ to $u_1,\cdots, u_d$ where $u_i\in V(P_i)$, $1\le i\le d$. It is possible that $x=u_i$ for $i\in[d]$, we may assume that $P_i=\{x\}$ under this circumstance. Analogously, there is a $d$-fan $Q_1, \cdots, Q_d$ in $G_n^2$ from $y$ to $v_1,\cdots, v_d$ where $v_i\in V(Q_i)$ and a $d$-fan $R_1, \cdots, R_d$ in $G_n^3$ from $z$ to $w_1,\cdots, w_d$ where $w_i\in V(R_i)$, $1\le i\le d$.

Note that $\{\hat{u}_i, \hat{v}_i, \hat{w}_i\}\subseteq V(G_n^{i+3})$ for $1\le i\le d$, there is a $\{\hat{u}_i, \hat{v}_i, \hat{w}_i\}$-tree $\widetilde{T}_i$ in $G_n^{i+3}$ since $G_n^{i+3}$ is connected.

For $1\le i\le d$, let
\begin{equation*}
T_{i}=\widetilde{T}_i\cup P_i\cup Q_i\cup R_i\cup\{u_i\hat{u}_i, v_i\hat{v}_i, w_i\hat{w}_i\}.
\end{equation*}
Then $T_1, \cdots,  T_{d}$ are $d$-internally edge disjoint $S$-trees in $H_n$. The proof is completed. \hfill$\Box$

\section{Applications} \label{secappli}

\subsection{Applications to the hierarchical star graph}\label{secHSn}

Let $\tau$ be a permutation on $[n]$, denote $\tau(1,i)$ be the permutation obtained by interchanging the 1st element with the $i$th element of $\tau$, where $2\le i\le n$.

\begin{definition}\label{defSn}(\cite{Sn1987})
An $n$-dimensional star graph, denoted by $S_n$, is an undirected graph with each vertex represented by a permutation on $[n]$ and two vertices $u$ and $v$ are adjacent if and only if $u=v(1,i)$ for some $i\in [n]\setminus \{1\}$.
\end{definition}

\begin{definition}\label{defHSn}(\cite{HSn2005})
For $n\ge 2$, a hierarchical star graph $HS_n$ of dimension $n$ consists of $n!$ $n$-dimensional star graphs $S_n$, called clusters. Each vertex of $HS_n$ is denoted by a two-tuple address $x=\langle c(x), p(x)\rangle$, where both $c(x)$ and $p(x)$ are arbitrary permutations on $[n]$. The first permutation $c(x)$ identifies the cluster the vertex $x$ belong to and the second permutation $p(x)$ identifies the vertex within the cluster. Two vertices $x=\langle c(x), p(x)\rangle$ and $y=\langle c(y), p(y)\rangle$ are adjacent in $HS_n$ if and only if one of the following three conditions holds:\\
(1) $c(x)=c(y)$ and $p(x)=p(y)(1,i)$ for some $2\le i \le n$; \\
(2) $c(x)\ne c(y)$, $c(x)=p(x)$ and $c(y)=p(y)=c(x)(1,n)$; \\
(3) $c(x)\ne c(y)$, $c(x)\ne p(x)$, $c(x)=p(y)$ and $p(x)=c(y)$.
\end{definition}

\begin{figure}[htbp]
\begin{minipage}[t]{0.2\linewidth}
\centering
\resizebox{0.9\textwidth}{!} {\includegraphics{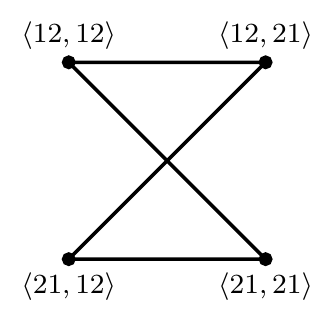}}
\caption{\small The hierarchical star graph $HS_2$} \label{figHS2}
\end{minipage}
\begin{minipage}[t]{0.8\linewidth}
\centering
\resizebox{0.9\textwidth}{!} {\includegraphics{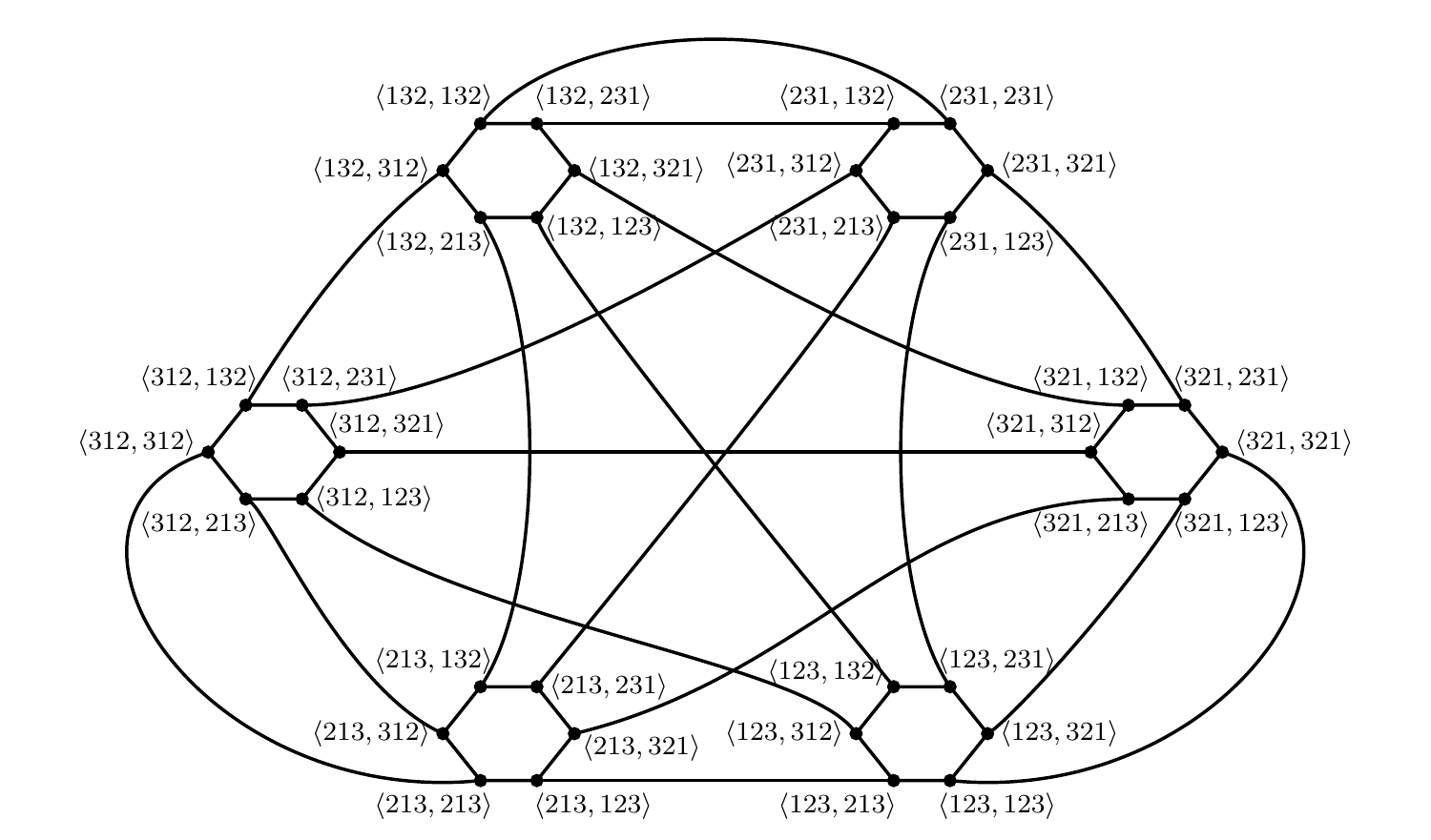}}
\caption{\small The hierarchical star graph $HS_3$} \label{figHS3}
\end{minipage}
\end{figure}

The hierarchical star graphs $HS_2$ and $HS_3$ are depicted in Figure \ref{figHS2} and Figure \ref{figHS3}, respectively. Note that the edges derived from the first condition of Definition \ref{defHSn} forms $n!$ vertex-disjoint subgraphs $S_n^i$ ($1\le i\le n!$), where each $S_n^i$ is isomorphic to the star graph $S_n$.

\begin{lemma}\label{lemSn}(\cite{Sn1987,SLi2016})
For any integer $n\ge 2$, $S_n$ is an ($n-1$)-regular graph with $n!$ vertices. Moreover, $\kappa(S_n)=n-1$ and $\kappa_3(S_n)=n-2$.
\end{lemma}

\begin{lemma}\label{lemHSn}(\cite{HSn2005,GuHao2018})
For any integer $n\ge 2$, $HS_n$ is an $n$-regular $n$-connected graph, and there is one or two cross edges between any pair of clusters.
\end{lemma}

\begin{corollary}\label{corHSn}
For $n\ge 2$, $\kappa_3(HS_n)=n-1$.
\end{corollary}

\noindent{\bf Proof}\;\; For $n\ge 3$, we have $n!\ge n+2$. Therefore, $HS_n$ is a special kind of graph $H_n$ defined in Definition \ref{defHGn}. By Theorem \ref{thmHGn2}, Lemma \ref{lemSn} and Lemma \ref{lemHSn}, it follows that $\kappa_3(HS_n)=n-1$ for $n\ge 3$.

Next we only need to show that $\kappa_3(HS_2)=1$. Firstly, it has $\kappa_3(HS_2)\le 1$ by Lemma \ref{lemupperK} and the fact that $HS_2$ is 2-regular. Secondly, it is easily seen that $\kappa_3(HS_2)\ge 1$ since $HS_2$ is connected.  The proof is completed. \hfill$\Box$

\subsection{Applications to the hierarchical cubic network}\label{secHCNn}

For any integer $n\ge 2$, the {\it $n$-dimensional hypercube}, denoted by $Q_n$, is the graph in which each vertex $x$ is corresponding to a distinct $n$-digit binary string $x_1x_2\cdots x_n$ on the set $\{0,1\}$, and two vertices $x$ and $y$ are adjacent in $Q_n$ if and only if $d_H(x,y)=1$, where $d_H(x,y)$ is the Hamming distance between $x$ and $y$. Let $x=x_1x_2\cdots x_n$ be an $n$-digit binary string, denote $\overline{x}=\overline{x_1} ~\overline{x_2}\cdots \overline{x_n}$, where $\overline{x_i}=1-x_i$ for all $i\in [n]$.

\begin{definition}\label{defHCNn}(\cite{HCNn1995})
The {\it $n$-dimensional hierarchical cubic network} $HCN_n$ ($n\ge 2$) can be decomposed into $2^n$ clusters, say $C_1, C_2, \cdots, C_{2^n}$, each cluster is isomorphic to an $n$-dimensional hypercube $Q_n$. Each vertex $x$ of $HCN_n$ is denoted by a two-tuple address $x=\langle c(x), p(x)\rangle$, where both $c(x)$ and $p(x)$ are $n$-digit binary strings. The first $n$-digit binary string $c(x)$ identifies the cluster the vertex $x$ belong to and the second $n$-digit binary string $p(x)$ identifies the vertex within the cluster. Two vertices $x=\langle c(x), p(x)\rangle$ and $y=\langle c(y), p(y)\rangle$ are adjacent in $HCN_n$ if and only if one of the following three conditions holds:\\
(1) $c(x)=c(y)$ and $d_H(p(x),p(y))=1$; \\
(2) $c(x)\ne c(y)$, $c(x)=p(x)$ and $c(y)=p(y)=\overline{c(x)}$; \\
(3) $c(x)\ne c(y)$, $c(x)\ne p(x)$, $c(x)=p(y)$ and $p(x)=c(y)$.
\end{definition}

\begin{lemma}\label{lemQn}(\cite{Bondy,HLiLi2012,K4Qn2017})
For any integer $n\ge 2$, the hypercube $Q_n$ is an $n$-regular graph with $2^n$ vertices. Furthermore, $\kappa(Q_n)=n$ and $\kappa_3(Q_n)=n-1$.
\end{lemma}

\begin{lemma}\label{lemHCN1}(\cite{HCNn1995,Cheng2014HCN,K4HCN2021})
The hierarchical cubic network $HCN_n$ ($n\ge 2$) has the following properties:\\
(1) $HCN_n$ is ($n+1$)-regular and ($n+1$)-connected;\\
(2) there is one or two cross edges between different clusters $C_i$ and $C_j$, ($i,j\in [2^n]$).
\end{lemma}

\begin{figure}[htbp]
\begin{minipage}[t]{0.5\linewidth}
\centering
\resizebox{0.8\textwidth}{!} {\includegraphics{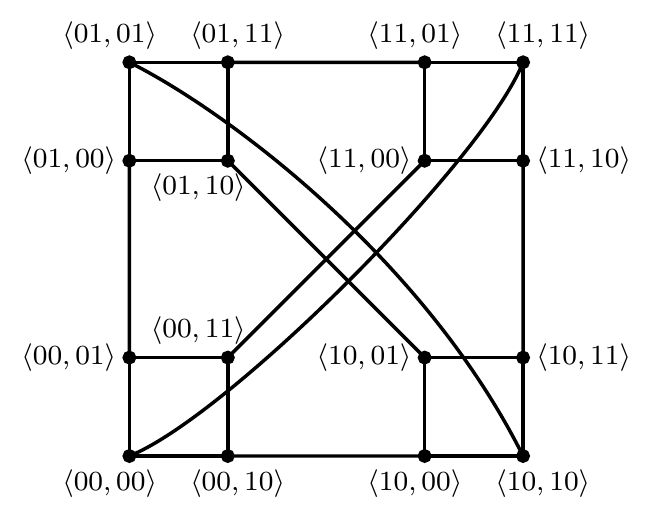}}
\caption{\small The hierarchical cubic network $HCN_2$} \label{figHCN2}
\end{minipage}
\quad
\begin{minipage}[t]{0.5\linewidth}
\centering
\resizebox{0.8\textwidth}{!} {\includegraphics{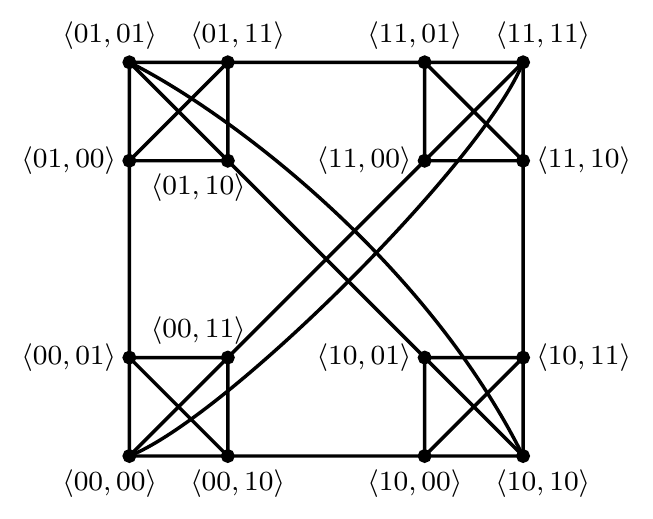}}
\caption{\small The hierarchical folded hypercube $HFQ_2$} \label{figHFQ2}
\end{minipage}
\end{figure}

The hierarchical cubic network $HCN_2$ is depicted in Figure \ref{figHCN2}.

\begin{lemma}\label{lemHCN2}
$\kappa_3(HCN_2)=2$.
\end{lemma}

\noindent{\bf Proof}\;\; Firstly, Lemma \ref{lemupperK} together with Figure \ref{figHCN2} enforce that $\kappa_3(HCN_2)\le \delta(HCN_2)-1=2$. By Lemma \ref{lemHCN1}, $\kappa(HCN_2)=3$. Therefore, $\kappa_3(HCN_2)\ge \lceil\frac{3}{2}\rceil=2$ according to Lemma \ref{lemkr}. \hfill $\Box$

\begin{corollary}\label{corHCNn}
For $n\ge 2$, $\kappa_3(HCN_n)=n$.
\end{corollary}

\noindent{\bf Proof}\;\; Since $2^n\ge n+3$ for $n\ge 3$. By Theorem \ref{thmHGn2}, Lemma \ref{lemQn} and Lemma \ref{lemHCN1}, it has $\kappa_3(HCN_n)=n$ for $n\ge 3$ (\cite{K4HCN2021}). Combined with Lemma \ref{lemHCN2}, it has $\kappa_3(HCN_n)=n$ for $n\ge 2$. \hfill$\Box$

\subsection{Applications to the hierarchical folded hypercube graph}\label{secHFQn}

For $n\ge 2$, the {\it $n$-dimensional folded hypercube} $FQ_n$ is a graph obtained from the hypercube $Q_n$ by adding an edge between any two vertices $x$ and $\overline{x}$ (\cite{El1991}).

\begin{definition}\label{defHFQn}(\cite{HFQn2000})
The {\it $n$-dimensional hierarchical folded cube} $HFQ_n$ ($n\ge 2$) can be decomposed into $2^n$ clusters, say $C_1, C_2, \cdots, C_{2^n}$, each cluster is isomorphic to an $n$-dimensional folded hypercube $FQ_n$. Each vertex $x$ of $HFQ_n$ is denoted by a two-tuple address $x=\langle c(x), p(x)\rangle$, where both $c(x)$ and $p(x)$ are $n$-digit binary strings. The first binary string $c(x)$ identifies the cluster the vertex $x$ belong to and the second binary string $p(x)$ identifies the vertex within the cluster. Two vertices $x=\langle c(x), p(x)\rangle$ and $y=\langle c(y), p(y)\rangle$ are adjacent in $HFQ_n$ if and only if one of the following four conditions holds:\\
(1) $c(x)=c(y)$ and $d_H(p(x),p(y))=1$; \\
(2) $c(x)=c(y)$ and $p(x)=\overline{p(y)}$; \\
(3) $c(x)\ne c(y)$, $c(x)=p(x)$ and $c(y)=p(y)=\overline{c(x)}$; \\
(4) $c(x)\ne c(y)$, $c(x)\ne p(x)$, $c(x)=p(y)$ and $p(x)=c(y)$.
\end{definition}

\begin{lemma}\label{lemFQn}(\cite{El1991,WangLi2022})
For $n\ge 2$, the folded hypercube $FQ_n$ is an ($n+1$)-regular graph with $2^n$ vertices. Moreover, $\kappa(FQ_n)=n+1$, $\kappa_3(FQ_n)=n$.
\end{lemma}

\begin{lemma}\label{lemHFQn}(\cite{HFQn2000,SunKHFQn2019})
The hierarchical folded hypercube $HFQ_n$ ($n\ge 2$) has the following properties:\\
(1) $HFQ_n$ is ($n+2$)-regular;\\
(2) there is one or two cross edges between different clusters $C_i$ and $C_j$, ($i,j\in [2^n]$);\\
(3) $\kappa(HFQ_n)=n+2$ for $n\ge 3$.
\end{lemma}

The hierarchical folded hypercube $HFQ_2$ is depicted in Figure \ref{figHFQ2}.

\begin{lemma}\label{remHFQ2}
$\kappa(HFQ_2)=4$ and $\kappa_3(HFQ_2)=3$.
\end{lemma}

\noindent {\bf Proof}\;\; First of all, by using almost the same arguments to that of Lemma 2.4 in \cite{SunKHFQn2019}, we can get that $\kappa(HFQ_2)=4$.

Now we shall prove that $\kappa_3(HFQ_2)=3$. Lemma \ref{lemupperK} and Lemma \ref{lemHFQn} imply that $\kappa_3(HFQ_2)\le \delta(HFQ_2)-1=3$. Moreover, $\kappa(HFQ_2)=4$ yields that $\kappa_3(HFQ_2)\ge 3$ according to Lemma \ref{lemkr}. \hfill $\Box$

\begin{corollary}\label{corHFQn}
For $n\ge 2$, $\kappa_3(HFQ_n)=n+1$.
\end{corollary}

\noindent{\bf Proof}\;\; Since $2^n\ge n+4$ for $n\ge 3$.  According to Theorem \ref{thmHGn2}, Lemma \ref{lemFQn} and Lemma \ref{lemHFQn}, we have $\kappa_3(HFQ_n)=n+1$ for $n\ge 3$. Combined with Lemma \ref{remHFQ2}, the result holds. \hfill$\Box$

\section{Conclusion}\label{secconclusion}

The generalized $k$-connectivity is a natural generalization of the traditional connectivity and can serve for measuring the capability of a network $G$ to connect any $k$ vertices in $G$. In this paper, we firstly introduce a family of regular networks and determine their generalized 3-connectivity. As applications, the generalized 3-connectivity of the hierarchical star graph $HS_n$, the hierarchical cubic network $HCN_n$ and the hierarchical folded hypercube $HFQ_n$, are determined.  We can see that most of the results on the generalized $k$-connectivity of networks are about $k=3$. It would be an interesting and challenging topic to study the generalized $k$-connectivity of $HS_n$ and $HFQ_n$ for $k\ge 4$.

%


\begin{thebibliography}{20}
%
%




\bibitem{Whitney1932} H.Whitney, Congruent graphs and connectivity of graphs, J. Amer. Math. Soc., 54 (1932) 150-168

\bibitem{Chartrand1984} G.Chartrand, S.F.Kapoor, L.Lesniak, D.R.Lick, Generalized connectivity in graphs, Bombay Math., 2 (1984) 1-6

\bibitem{SLi2012n} S.Li, X.Li, Note on the hardness of generalized connectivity, J. Comb. Optim., 24 (2012) 389-396

\bibitem{HZLi2014} H.Li, X.Li, Y.Mao, Y.Sun, Note on the generalized connectivity, Ars Combin., 114 (2014) 193-202

\bibitem{SLi2010} S.Li, X.Li, W.Zhou, Sharp bounds for the generalized connectivity $\kappa_3(G)$, Discrete Math., 310 (2010) 2147-2163

\bibitem{Chartrand2010} G. Chartrand, F. Okamoto, P. Zhang, Rainbow trees in graphs and generalized connectivity, Networks, 55 (4) (2010) 360-367

\bibitem{SLi2012b} S.Li, W.Li, X.Li, The generalized connectivity of complete bipartite graphs, Ars Combin., 104 (2012) 65-79

\bibitem{HZLi2012} H.Li, X.Li, Y.Sun, The generalized 3-connectivity of Cartesian product graphs, Discrete Math. Theor. Comput. Sci., 14 (1) (2012) 43-54

\bibitem{HZLi2017} H.Li, Y.Ma, W.Yang, Y.Wang, The generalized 3-connectivity of graph products, Appl. Math. Comput., 295 (2017) 77-83


\bibitem{SLin2017} S.Lin, Q.Zhang, The generalized 4-connectivity of hypercubes, Discrete Appl. Math., 220 (2017) 60-67

\bibitem{ZhaoHao20191} S.Zhao, R.Hao, E.Cheng, Two kinds of generalized connectivity of dual cubes, Discrete Appl. Math., 257 (2019) 306-316

\bibitem{ZhaoHao20192} S.Zhao, R.Hao, The generalized 4-connectivity of exchanged hypercubes, Applied Math. Comput., 347 (2019) 342-353

\bibitem{Wei2021} C.Wei, R.Hao, J.Chang, The reliability analysis based on the generalized connectivity in balanced hypercubes, Discrete Appl. Math., 292 (2021) 19-32

\bibitem{Wang2021} J.Wang, The generalized 3-connectivity of two kinds of regular networks, Theoret. Comput. Sci., 893 (2021) 183-190

\bibitem{SLi2017} S.Li, Y.Shi, J.Tu, The generalized 3-connectivity of Cayley graphs on symmetric groups generated by trees and cycles, Graphs and Combin., 33 (2017) 1195-1209

\bibitem{ZHao20193} S.Zhao, R.Hao, The generalized three-connectivity of two kinds of Cayley graphs, The Comput. Journal, 62 (2019) 144-149

\bibitem{SLi2016} S.Li, J.Tu, C.Yu, The generalized 3-connectivity of star graphs and bubble-sort graphs, Appl. Math. Comput., 274 (2016) 41-46

\bibitem{Hao20191} S.Zhao, R.Hao, The generalized connectivity of bubble-sort star graphs, International J. Foundations of Comput. Sci., 30 (2019) 793-809

\bibitem{HSn2005} W.Shi, P.K.Srimani, Hierachical star: a new two level interconnection network, J. Syst. Archit., 51 (2005) 1-14

\bibitem{HCNn1995} K.Ghose, K.R.Desai, Hierarchical cubic network, IEEE Trans. Parallel Distrib. Syst., 6 (4) (1995) 427-435

\bibitem{HFQn2000} Y.Shi, Z.Hou, J.Song, Hierarchical interconnection networks with folded hypercubes as basic cluster, in: Proceedings of the 4th International Conference, Exhibition on High Performance Computing in the Asia-Pacific Region 1, (2000) 134-137

\bibitem{Bondy} J.A.Bondy, U.S.R.Murty, Graph theory, Springer, New York, 2007

\bibitem{Sn1987} S.B.Akers, B.Krishnamurthy, D.Harel, The star graph: an attractive alternative to the $n$-cube, in: Proc. Int'l Conf. Parallel Process, 1987, 393-400

\bibitem{GuHao2018} M.Gu, R.Hao, L.Jiang, Fault-tolerance and diagnosability of hierarchical star networks, International J. Comput. Math.: Comput. Sys. Theory, 3 (2018) 106-121

\bibitem{HLiLi2012} H.Li, X.Li, Y.Sun, The generalized 3-connectivity of Cartesian product graphs, Discret. Math. Theor. Comput. Sci., 14(1) (2012) 43-54

\bibitem{K4Qn2017} S.Lin, Q.Zhang, The generalized 4-connectivity of hypercubes, Discrete Appl. Math., 220 (2017) 60-67

\bibitem{Cheng2014HCN} E.Cheng, K.Qiu, Z.Shen, Connectivity results of hierarchical cubic networks as associated with linearly many faults, in: IEEE 17th International Conference on Computational Science and Engineering, (2014) 1213-1220


\bibitem{K4HCN2021} S.Zhao, R.Hao, J.Wu, The generalized 4-connectivity of hierarchical cubic networks, Discrete Appl. Math., 289 (2021) 194-206

\bibitem{El1991} A.El-Amawy, S.Latifi, Properties and performance of folded hypercubes, IEEE Trans. Parallel Distrib. Syst., 2 (1991) 31-42

\bibitem{WangLi2022} J.Wang, F.Li, The generalized 3-connectivity of the folded hypercube $FQ_n$, The Comput. Journal, https://doi.org/10.1093/comjnl/bxac137

\bibitem{SunKHFQn2019} X.Sun, Q.Dong, S.Zhou, et al, Fault tolerance analysis of hierarchical folded cube, Theore. Comput. Sci., 790 (2019) 117-130




\end{thebibliography}


\end{document}